\documentclass[12pt, a4paper,revtex4]{article}
\usepackage{amssymb, amsfonts, amsmath, latexsym, amsthm, yfonts}

\newtheorem{question}{Question}
\theoremstyle{definition}
\newtheorem{theorem}{}[section]

\begin{document}
\date{}
\title{\bf Combinatorial Derivation}
\author{Igor~V.~Protasov}

\maketitle

\begin{abstract}
Let $G$ be a group, $\mathcal{P}_G$ be the family of all subsets of $G$. For a subset $A\subseteq G$, we put
$\Delta(A)=\{g\in G:|gA\cap A|=\infty\}$. The mapping $\Delta:\mathcal{P}_G\rightarrow\mathcal{P}_G$, $A\mapsto\Delta(A)$, is called a combinatorial derivation and can be considered as an analogue of the topological derivation $d:\mathcal{P}_X\rightarrow\mathcal{P}_X$, $A\mapsto A^d$, where $X$ is a topological space and $A^d$ is the set of all limit points of $A$. Content: elementary properties, thin and almost thin subsets, partitions, inverse construction and $\Delta$-trajectories,  $\Delta$ and $d$.

\

\textbf{Classification}. 20A05, 20F99, 22A15, 06E15, 06E25.

\


\end{abstract}

Let $G$ be a group with the identity $e$, $\mathcal{P}_G$ be the family of all subsets of $G$. For a subset $A$ of $G$, we denote
$$\Delta(A)=\{g\in G:|gA\cap A|=\infty\},$$
observe that $\Delta(A)\subseteq AA^{-1}$, and say that the mapping
$$\Delta:\mathcal{P}_G\rightarrow\mathcal{P}_G, A\mapsto\Delta(A)$$
is a {\it combinatorial derivation}.

In this paper, on one hand, we analyze from the $\Delta$-point of view a series of results from Subset Combinatorics of Groups (see the survey \cite{b9}), and point out some directions for further progress. On the other hand, the $\Delta$-operation is interesting and intriguing by its own sake. In contrast to the trajectory 
$A\rightarrow AA^{-1}\rightarrow (AA^{-1})(AA^{-1})\rightarrow\ldots$, 
the $\Delta$-trajectory $A\rightarrow\Delta(A)\rightarrow\Delta^2(A)\rightarrow\ldots$ of a subset $A$ of $G$ could be surprisingly complicated: stabilizing, increasing, decreasing, periodic or chaotic. For a symmetric subset $A$ of $G$ with $e\in A$, there exists a subset $X\subseteq G$ such that $\Delta(X)=A$. We conclude the paper by demonstrating how $\Delta$ and a topological derivation $d$ arise from some unified ultrafilter construction.

We note also that $\Delta(A)$ may be considered as some infinite version of the symmetry sets well-known in Additive Combinatorics \cite[p. 84]{b13}. Given a finite subset $A$ of an Abelian group $G$ and $\alpha \geqslant 0$, the symmetry set $Sym_\alpha(A)$ is defined by

$$Sym_\alpha(A) = \{g\in G: |A\cap (A+g)| \geqslant \alpha|A|\}.$$

\section{Elementary properties}{\label{s1}}

\begin{theorem}{\label{t11}}
$(\Delta(A))^{-1}=\Delta(A)$, $\Delta(A)\subseteq AA^{-1}$.
\end{theorem}

\begin{theorem}{\label{t12}}
$\Delta(A)=\varnothing\Leftrightarrow e\notin\Delta(A)\Leftrightarrow A\text{ is finite}$.
\end{theorem}

\begin{theorem}{\label{t13}}
For subsets $A,B$ of $G$, we let
$$\Delta(A,B)=\{g\in G:|gA\cap B|=\infty\}$$
and note that
$$\Delta(A\cup B)=\Delta(A)\cup\Delta(B)\cup\Delta(A,B)\cup\Delta(B,A),$$
$$\Delta(A\cap B)\subseteq\Delta(A)\cap\Delta(B)$$.
\end{theorem}

\begin{theorem}{\label{t14}}
If $F$ is a finite subset of $G$ then
$$\Delta(FA)=F\Delta(A)F^{-1}.$$
\end{theorem}

\begin{theorem}{\label{t15}}
If $A$ is an infinite subgroup then $A=\Delta(A)$ but the converse statement does not hold, see \textbf{\ref{t54}}.
\end{theorem}

\section{Thin and almost thin subsets}{\label{s2}}

A subset $A$ of a group $G$ is said to be \cite{b8}:
\begin{itemize}
\item \textit{thin} if either $A$ is finite or $\Delta(A)=\{e\}$;
\item \textit{almost thin} if $\Delta(A)$ is finite;
\item \textit{$k$-thin} ($k\in\mathbb{N}$) if $|gA\cap A|\leqslant k$ for each $g\in G\setminus\{e\}$;
\item \textit{sparse} if, for every infinite subset $X\subseteq G$, there exists a non-empty finite subset $F\subset X$ such that $\bigcap_{g\in F}gA$ is infinite;
\item \textit{$k$-sparse} ($k\in\mathbb{N}$) if, for every infinite subset $X\subseteq G$, there exists a subset $F\subset X$ such that $|F|\leqslant k$ and $\bigcap_{g\in F}gA$ is finite.
\end{itemize}

The following statements are from \cite{b8}.

\begin{theorem}{\label{t21}}
Every almost thin subset $A$ of a group $G$ can be partitioned in $3^{|\Delta(A)|-1}$ thin subsets. If $G$ has no elements of odd order, then $A$ can be partitioned in $2^{|\Delta{A}|-1}$ thin subsets.
\end{theorem}

\begin{theorem}{\label{t22}}
A subset $A$ of a group $G$ is 2-sparse if and only if $X^{-1}X\nsubseteq\Delta(A)$ for every infinite subset $X$ of $G$.
\end{theorem}

\begin{theorem}{\label{t23}}
For every countable thin subset $A$ of a group $G$, there is a thin subset $B$ such that $A\cup B$ is 2-sparse but not almost thin.
\end{theorem}

\begin{theorem}{\label{t24}}
Suppose that a group $G$ is either torsion-free or, for every $n\in\mathbb{N}$, there exists a finite subgroup $H_n$ of $G$ such that $|H_n|>n$. Then there exists a 2-sparse subset of $G$ which cannot be partitioned in finitely many thin subsets.
\end{theorem}

By \textbf{\ref{t22}}, every almost thin subset is 2-sparse. By \textbf{\ref{t23}}, \textbf{\ref{t24}}, the class of 2-sparse subsets is wider than the class of almost thin subsets. By \textbf{\ref{t23}}, a union of two thin subsets needs not to be almost thin. By \textbf{\ref{t13}}, a union $A_1\cup\ldots\cup A_n$ of almost thin subset is almost thin if and only if $\Delta(A_i,A_j)$ is finite for all $i,j\in\{1,\ldots,n\}$, By \textbf{\ref{t14}}, if $A$ is almost thin and $K$ is finite then $KA$ is almost thin.

The following statements are from \cite{b7}.

\begin{theorem}{\label{t25}}
For every infinite group $G$, there exists a 2-thin subset such that $G=XX^{-1}\cup X^{-1}X$.
\end{theorem}

\begin{theorem}{\label{t26}}
For every infinite group $G$, there exists a 4-thin subset such that $G=XX^{-1}$.
\end{theorem}

Since $\Delta(X)=\{e\}$ for each infinite thin subset of $G$, \textbf{\ref{t26}} gives us $X$ with $\Delta(X)=\{e\}$ and $XX^{-1}=G$.

\section{Large and small subsets}{\label{s3}}

A subset $A$ of a group $G$ is called \cite{b8}:
\begin{itemize}
\item \textit{large} if there exists a finite subset $F$ of $G$ such that $G=FA$;
\item \textit{$\Delta$-large} if $\Delta(A)$ is large;
\item \textit{small} if $(G\setminus A)\cap L$ is large for each large subset $L$ of $G$;
\item \textit{P-small} if there exists an injective sequence $(g_n)_{n\in\omega}$ in $G$ such that the subsets $\{g_nA: n\in\omega\}$ are pairwise disjoint;
\item \textit{almost P-small} if there exists an injective sequence $(g_n)_{n\in\omega}$ in $G$ such that the family $\{g_nA: n\in\omega\}$ is almost disjoint, i.e. $g_nA\cap g_mA$ is finite for all distinct $n,m\in\omega$.
\item \textit{weakly P-small} if, for every $n\in\omega$, one can find distinct elements $g_1,\ldots,g_n$ of $G$ such that the subsets $g_1A,\ldots,g_nA$ are pairwise disjoint.
\end{itemize}

\begin{theorem}{\label{t31}}
Let $G$ be a group, $A$ is a large subset of $G$. We take a finite subset $F$ of $G$, $F=\{g_1,\ldots, g_n\}$ such that $G=FA$. Take an arbitrary $g\in G$. Then $g_iA\cap gA$ is infinite for some $i\in \{1,\ldots, n\}$, so $g_i^{-1}g\in \Delta(A)$. Hence, $G=F\Delta(A)$ and $A$ is $\Delta$-large. By \textbf{\ref{t26}}, the converse statement is very far from being true.

If $A$ is not small then $FA$ is thick (see \textbf{\ref{t43}}) for some finite subset $F$. It follows that $\Delta(FA) = G$. By \textbf{\ref{t14}}, $\Delta(FA) = F\Delta(A) F^{-1}$, so if $G$ is Abelian then $A$ is $\Delta$-large.
\end{theorem}

\begin{question}
Is every nonsmall subset of an arbitrary infinite non-Abelian group $G$ $\Delta$-large?
\end{question}

\begin{theorem}{\label{t32}}
It is easy to see that $A$ is P-small (almost P-small) if and only if there exist an infinite subset $X$ of $G$ such that $X^{-1}X\cap PP^{-1}=\{e\}$ ($X^{-1}X\cap\Delta(X)=\{e\}$). $A$ is weakly P-small if and only if, for every $n\in\omega$, there exists $F\subset G$ such that $|F|=n$ and $F^{-1}F\cap PP^{-1}=\{e\}$
\end{theorem}

\begin{theorem}{\label{t33}}
By \cite[Lemma 4.2]{b8}, if $AA^{-1}$ is not large then $A$ is small and P-small. Using the inverse construction \textbf{\ref{t53}}, we can find $A$ such that $A$ is not $\Delta$-large and $A$ is not P-small.
\end{theorem}

\begin{theorem}{\label{t34}}
Every infinite group $G$ has a weakly P-small not P-small subsets \cite{b1}. Moreover, $G$ has almost P-small not P-small subset and , if $G$ is countable, weakly P-small not almost P-small subset. By \cite{b8}, every almost P-small subset can be partitioned in two P-small subsets. If $A$ is either almost or weakly P-small then $G\setminus\Delta(A)$ is infinite, but a subset $A$ with infinite $G\setminus\Delta(A)$ could be large: $G=\mathbb{Z}$, $A=2\mathbb{Z}$.
\end{theorem}

\section{Partitions}{\label{s4}}

\begin{theorem}{\label{t41}}
Let $G$ be a group and let $G=A_1\cup\ldots A_n$ be a finite partition of $G$. In section \ref{s6}, we show that at least one cell $A_i$ is $\Delta$-large, in particular, $A_iA_i^{-1}$ is large. If $G$ is infinite amenable group and $\mu$ is a left invariant Banach measure on $G$, we can strengthened this statement: there exist a cell $A_i$ and a finite subset $F$ such that $|F|\leqslant n$ and $G=F\Delta(A_i)$. To verify this statement, we take $A_i$ such that $\mu(A_i)\geqslant\frac{1}{n}$ and choose distinct $g_1,\ldots,g_m$ such that $\mu(g_kA_i\cap g_l A_i)=0$ for all distinct $k,l\in\{1,\ldots,m\}$, and the family $\{g_1A_i,\ldots,g_mA_i\}$ is maximal with respect to this property. Clearly, $m\leqslant n$. For each $g\in G$, we have $\mu(gA_i\cap g_kA_i)>0$ for some $k\in\{1,\ldots,m\}$ so $g_k^{-1}g\in\Delta(A_i)$ and $G=\{g_1,\ldots,g_m\}\Delta(A_i)$.
\end{theorem}

\begin{theorem}{\label{t42}}
By \cite[Theorem 12.7]{b10}, for every partition $A_1\cup\ldots\cup A_n$ of an arbitrary group $G$, there exist a cell $A_i$ and a finite subset $F$ of $G$ such that $G=FA_iA_i^{-1}$ and $|F|\leqslant 2^{2^{n-1}-1}$. S.~Slobodianiuk strengthened this statement: there are $F$ and $A_i$ such that $|F|\leqslant 2^{2^{n-1}-1}$ and $G=F\Delta(A_i)$.
\end{theorem}

It is an old unsolved problem \cite[Problem 13.44]{b5} whether $i$ and $F$ can be chosen so that $G=FA_iA_i^{-1}$ and $|F|\leqslant n$, see also \cite[Question 12.1]{b10}.

\begin{question}
Given any partition $G = A_1 \cup \ldots \cup A_n$, do there exist $F$ and $A_i$ such that $G=F \Delta(A_i)$ and $|F|\leqslant 2^n$?
\end{question}

\begin{theorem}{\label{t43}}
A subset $A$ of a group $G$ is called \cite{b11}:
\begin{itemize}
\item \textit{thick} if $G\setminus A$ is not large;
\item \textit{$k$-prethick} ($k\in\mathbb{N}$) if there exists a subset $F$ of $G$ such that $|F|\leqslant k$ and $FA$ is thick;
\item \textit{prethick} if $A$ is $k$-prethick for some $k\in\mathbb{N}$.
\end{itemize}
By \cite[Theorem 11.2]{b10}, if a group $G$ is partitioned into finite number of cells, then at least one cell of the partition is prethick. We say that a partition $\mathcal{P}$ of a group $G$ is {\it $k$-meager} if each number of $\mathcal{P}$ is not $k$-prethick. If $G$ is either countable locally finite, or countable residually finite, or infinite Abelian, then $G$ admits a $k$-meager 2-partition for each $k\in\mathbb{N}$ \cite{b11}.
\end{theorem}

\begin{theorem}{\label{t44}}
By \cite[Theorem 5.3.2]{b3}, for a group $G$, the following two conditions (i) and (ii) are equivalent:
\begin{itemize}
\item[(i)] for every partition $G=A\cup B$, either $G=AA^{-1}$ or $G=BB^{-1}$;
\item[(ii)] each element of $G$ has odd order.
\end{itemize}
If $G$ is infinite, we can show that these conditions are equivalent to
\begin{itemize}
\item[(iii)] for every partition $G=A\cup B$, either $G=\Delta(A)$ or $G=\Delta(G)$.
\end{itemize}
\end{theorem}

\section{Inverse construction and $\Delta$-trajectories}{\label{s5}}

\begin{theorem}{\label{t51}}
Let $G$ be an infinite group, $A\subseteq G$, $A=A^{-1}$, $e\in A$. We construct a subset $X$ of $G$ such that $\Delta(X)=A$.
\end{theorem}

First, assume that $G$ is countable and write the elements of $A$ in the list $\{a_n:n<\omega\}$, if $A$ is finite then all but finitely many $a_n$ are equal to $e$. We represent $G\setminus A$ as a union $G\setminus A=\bigcup_{n\in\omega}B_n$ of finite subsets such that $B_n\subseteq B_{n+1}$, $B_n^{-1}=B_n$. Then we choose inductively a sequence $(X_n)_{n\in\omega}$ of finite subsets of $G$,
$$X_n=\{x_{n0},x_{n1},\ldots,x_{nn}, a_0x_{n0}, \ldots,a_nx_{nn}\}$$
such that $X_mX_n^{-1}\cap B_n=\{e\}$ for all $m\leqslant n<\infty$.

After $\omega$ steps, we put $X=\bigcup_{n\in\omega}X_n$. By the construction, $\Delta(X)=A$.

If $|A|\leqslant\aleph_0$ but $G$ is not countable, we take a countable subgroup $H$ of $G$ such that $A\subseteq H$, forget about $G$ and find a subset $X\subseteq H$ such that $\Delta(X)$ is equal to $A$ in $H$. Since $gA\cap A=\varnothing$ for each $g\in G\setminus H$, we have $\Delta(X)=A$.

At last, let $|A|>\aleph_0$. By above paragraph, we may suppose that $|A|=|G|$. We enumerate $A=\{a_\alpha:\alpha<|G|\}$ and construct inductively a sequence $(X_\alpha)_{\alpha<|G|}$ of finite subsets of $G$ and an increasing sequence $(H_\alpha)_{\alpha<|G|}$ of subgroup of $G$ such that if $\alpha=0$ or $\alpha$ is a limit ordinal, $n\in\omega$,
$$X_{\alpha+n}=\{x_{\alpha+n,0},x_{\alpha+n,1},\ldots,x_{\alpha+n,n},a_\alpha x_{\alpha+n,0},\ldots, a_{\alpha+n}x_{\alpha+n,n}\},$$
$$X_{\alpha+n}\subseteq H_{\alpha+n+1}\setminus H_{\alpha+n},\ X_{\alpha+n}X_{\alpha+n}^{-1}\subseteq A\cup(H_{\alpha+n+1}\setminus H_{\alpha+n}).$$
After $|G|$ steps, we put $X=\bigcup_{\alpha<|G|}X_\alpha$. By the construction, $\Delta(X)=A$.

\begin{theorem}{\label{t52}}
Let $A_1,\ldots,A_m$ be subsets of an infinite group $G$ such that $G=A_1\cup\ldots\cup A_m$. By the Hindman theorem \cite[Theorem 5.8]{b4}, there are exists $i\in\{1,\ldots,m\}$ and an injective sequence $(g_n)_{n\in\omega}$ in $G$ such that $FP(g_n)_{n\in\omega}\subseteq A_i$, where $FP(g_n)_{n\in\omega}$ is a set of all element of the form $g_{i_1}g_{i_2}\ldots g_{i_l}$, $i_1<\ldots,i_k<\omega$, $k\in\omega$.
\end{theorem}

We show that there exists $X\subseteq FP(g_n)_{n\in \omega}$ such that $\Delta(X) = \{e\} \cup FP(g_n)_{n\in \omega} \cup (FP(g_n)_{n\in \omega})^{-1}$. We note that if $G$ is countable, at each step $n$ of the inverse construction \ref{t52}, the elements $x_{n0}, \ldots , x_{nn}$ can be chosen from any pregiven infinite subset $Y$ of $G$. We enumerate $FP(g_n)_{n\in \omega}$ in a sequence $(a_n)_{n\in \omega}$ and put $Y=\{g_n: n\in \omega\}$. Using above observation, we get the desired $X$.

\begin{theorem}{\label{t53}}
If $G$ is countable, we can modify the inverse construction to get $X$ such that $\Delta(X)=A$ and $|X\cap g_1\cap g_2X|<\infty$ for all distinct $g_1,g_2\in G\setminus\{e\}$, in particular, $X$ is 3-sparse and, in particular, small.
\end{theorem}

Another modification, we can choose $X$ such that $X\cap gX\ne\varnothing$ for each $g\in G$. If we take $A$ not large, then we get $X$ which is not P-small and $X$ is not $\Delta$-large, see \textbf{\ref{t33}}.

\begin{theorem}{\label{t54}}
Let $G$ be a countable group such that, for each $g\in G\setminus\{e\}$, the set $\sqrt{g}=\{x\in G:x^2=g\}$ is finite. The following statements show that all possible $\Delta$-trajectories of subsets of $G$ can be realized.
\begin{itemize}
\item[($Tr_1$)] Given any subset $X_0\subseteq G$, $X_0=X_0^{-1}$, $e\in X_0$, there exists a sequence $(X_n)_{n\in\omega}$ of subsets of $G$ such that $\Delta(X_{n+1})=X_n$ and $X_m\cap X_n=\{e\}$, $0 <m<n<\omega$.
\item[($Tr_2$)] There exists a sequence $(X_n)_{n\in\mathbb{Z}}$ of subsets of $G$ such that $\Delta(X_n)=X_{n+1}$, $X_m\cap X_n=\{e\}$, $m,n \in \mathbb{Z}, m\ne n$.
\item[($Tr_3$)] There exists a subset $A$ of $G$ such that $\Delta(A)=A$ but $A$ is not a subgroup.
\item[($Tr_4$)] There exists a subset $A$ such that $A\supset\Delta(A)\supset\Delta^2(A)\supset\ldots$.
\item[($Tr_5$)] There exists a subset $A$ such that $A\subset\Delta(A)\subset\Delta^2(A)\subset\ldots$.
\item[($Tr_6$)] For each natural natural number $n$, there exists a periodic $\Delta$-trajectory $X_0,\ldots,X_{n-1}$ of length $n$: $X_1=\Delta(X_0),X_2=\Delta(X_1),\ldots,X_n=\Delta(X_{n-1})$ such that $X_i\cap X_j=\{e\}$, $i<j<n$.
    
\end{itemize}
\end{theorem}

We use the following simple observation

\begin{itemize}
\item[(*)] if $F$ is a finite subset of an infinite group $G$ and $g\notin F$ then the set $\{x\in G: x^{-1}gx \notin F\}$ is infinite.
\end{itemize}

In constructions of corresponding trajectories, at each inductive step, we use a finiteness of $\sqrt{g}$ and (*) in the following form:
\begin{itemize}
\item[(**)] if $a\in G$, $F$ is a finite subset of $G$, $F\cap\{e,a^{\pm 1}\}=\varnothing$ then there exists $x\in G$ such that
    $$\{x^{\pm 1},(ax)^{\pm 1}\}\{x^{\pm 1}, (ax)^{\pm 1}\}\cap F = \varnothing.$$
\end{itemize}

We show how to get a 2-periodic trajectory: $X$, $Y$, $\Delta(X)=Y$, $\Delta(Y)=X$, $X\cap Y=\{e\}$. We write $G$ as a union $G=\bigcup_{n\in\omega}F_n$  of increasing chain $\{F_n:n\in\omega\}$ of finite symmetric subsets $F_0=\{e\}$. We put $X_0=Y_0=\{e\}$ and construct inductively with usage of \textit{(**)} two chains $(X_n)_{n\in\omega}$, $(Y_n)_{n\in\omega}$ of finite subsets of $G$ such that, for each $n\in\omega$,
$$X_{n+1}=\{(x(y))^{\pm 1},(yx(y))^{\pm 1}: y\in Y_0\cup\ldots\cup Y_n\},$$
$$Y_{n+1}=\{(y(x))^{\pm 1},(xy(x))^{\pm 1}: x\in X_0\cup\ldots\cup X_n\},$$
$$(X_0\cup\ldots\cup X_n)\cap(Y_0\cup\ldots\cup Y_n)=\{e\},$$
$$X_{n+1}X_{n+1}\cap(F_{n+1}\setminus(Y_0\cup\ldots\cup Y_n))=\varnothing,$$
$$Y_{n+1}Y_{n+1}\cap(F_{n+1}\setminus(X_0\cup\ldots\cup X_n))=\varnothing,$$
$$(X_0\cup\ldots\cup X_n)X_{n+1}\cap(F_{n+1}\setminus(Y_0\cup\ldots\cup Y_n))=\varnothing,$$
$$(Y_0\cup\ldots\cup Y_n)Y_{n+1}\cap(F_{n+1}\setminus(X_0\cup\ldots\cup X_n))=\varnothing.$$
After $\omega$ steps, we put $X=\bigcup_{n\in\omega}X_n$, $Y=\bigcup_{n\in\omega}Y_n$.

\section{$\Delta$ and $d$}{\label{s6}}

For a subset $A$ of a topological space $X$, the subset $A^d$ of all limit points of $A$ is called \textit{a derived subset}, and the mapping $d:\mathcal{P}(X)\rightarrow \mathcal{P}(X)$, $A\rightarrow A^d$, defined on the family of $\mathcal{P}(X)$ of all subsets of $X$, is called \textit{a topological derivation}, see \cite[\S9]{b6}.

Let $X$ be a discrete set, $\beta X$ be the Stone-\v{C}ech compactification of $X$. We identify $\beta X$ with the set of all ultrafilters on $X$, $X$ with the set of all principal ultrafilters, and denote $X^*=\beta X\setminus X$ the set of all free ultrafilters. The topology of $\beta X$ can be defined by the family $\{\overline{A}:A\subseteq X\}$ as a base for open sets, $\overline{A}=\{p\in\beta X:A\in p\}$, $A^*=\overline{A}\cap G^*$. For a filter $\varphi$ on $X$, we put $\overline{\varphi}=\{p\in\beta X:\varphi\subseteq p\}$, $\varphi^*=\overline{\varphi}\cap G^*$.

Let $G$ be a discrete group, $p\in\beta G$. Following \cite[Chapter 3]{b2}, we denote
$$cl(A,p)=\{g\in G:A\in gp\},\ gp=\{gP:P\in p\},$$
say that $cl(A,p)$ is \textit{a closure of $A$ in the direction of $p$}, and note that
$$\Delta(A)=\bigcap_{p\in A^*}cl(A,p).$$

A topology $\tau$ on a group $G$ is called \textit{left invariant} if the mapping $l_g:G\rightarrow G$, $l_g(x)=gx$ is continuous for each $g\in G$. A group $G$ endowed with a left invariant topology $\tau$ is called \textit{left topological}.
We note that a left invariant topology $\tau$ on $G$ is uniquely determined by the filter $\varphi$ of neighbourhoods of the identity $e\in G$, $\overline{\varphi}$ and $\varphi^*$ are the sets of all ultrafilters an all free ultrafilters of $G$ converging to $e$. For a subset $A$ of $G$, we have
$$A^d=\bigcap_{p\in(\tau^*)}cl(A,p),$$
and note that $A^d\subseteq\Delta(A)$ if $A$ is a neighbourhood of $e$ in $(G,\tau)$.

Now we endow $G$ with the discrete topology and, following \cite[Chapter 4]{b4}, extend the multiplication on $G$ to $\beta G$. For $p,q\in\beta G$, we take $P\in p$ and, for each $g\in P$, pick some $Q_g\in q$. Then $\bigcup_{g\in p}gQ_p\in pq$ and each member of $pq$ contains a subset of this form. With this multiplication, $\beta G$ is a compact right topological semigroup. The product $pq$ can also be defined by the rule \cite[Chapter 3]{b2}:
$$A\subseteq G,\ A\in pq \Leftrightarrow cl(A,q)\in p.$$

If $(G,\tau)$ is left topological semigroup then $\overline{\tau}$ is a subsemigroup of $\beta G$. If an ultrafilter $p\in\overline{\tau}$ is taken from the minimal ideal $K(\overline{\tau})$ of $\overline{\tau}$, by \cite[Theorem 5.0.25]{b2}. there exists $P\in p$ and finite subset $F$ of $G$ such that $Fcl(P,p)$ is neighbourhood of $e$ in $\tau$. In particular, if $\tau$ indiscrete ($\tau=\{\varnothing,G\}$), $p\in K(\beta G)$) and $P\in p$ then $cl(P,p)$ is large. If $G$ is infinite, $p\in K(\beta G)$ is free, so $cl(P,p)\subseteq\Delta(P)$ and $P$ is $\Delta$-large. If a group $G$ is finitely partitioned $G=A_1\cup\ldots\cup A_n$, then some cell $A_i$ is a member of $p$, hence $A_i$ is $\Delta$-large.


\newpage

\begin{thebibliography}{}

\bibitem{b1}
T.~Banakh, N.~Lyaskovska, 
{\it Weakly P-small not P-small subsets in groups}, 
Intern. J. Algebra Computations, \textbf{18} (2008), 1--6.

\bibitem{b2}
M.~Filali, I.~Protasov,
{\it Ultrafilters and Topologies on Groups},
Math. Stud. Monorg. Ser., Vol.~\textbf{13}, VNTL Publishers, Lviv, 2010.

\bibitem{b3}
V.~Gavrylkiv,
{\it Algebraic-topological structure on superextensions},
Dissertation, Lviv, 2009.

\bibitem{b4}
N.~Hindman, D.~Strauss,
{\it Algebra in the Stone-\v{C}ech compactification},
Walter de Grueter, Berlin, New York, 1998.

\bibitem{b5}
{\it The Kourovka Notebook},
Novosibirsk, Institute of Math., 1995.

\bibitem{b6}
K.~Kuratowski,
{\it Topology},
Vol.~\textbf{1}, Academic Press, New York and London, PWN, Warszawa, 1969.

\bibitem{b7}
Ie.~Lutsenko,
{\it Thin systems of generators of groups},
Algebra and Discrete Math., \textbf{9} (2010), 108--114.

\bibitem{b8}
Ie.~Lutsenko, I.~V.~Protasov,
{\it Sparse, thin and other subsets of groups},
Intern. J. Algebra Computation, \textbf{19} (2009), 491--510.

\bibitem{b9}
I.~V.~Protasov,
{\it Selective survey on Subset Combinatorics of Groups},
J. Math. Sciences, \textbf{174} (2011), 486--514.

\bibitem{b10}
I.~Protasov, T.~Banakh.,
{\it Ball Structure and Colorings of Groups and Graphs},
Math. Stud. Monorg. Ser., Vol.~\textbf{11}, VNTL Publishers, Lviv, 2003.

\bibitem{b11}
I.~Protasov, S.~Slobodianiuk,
{\it Prethick subsets in partitions of groups}, Centr. Europ. Math. J., submitted.

\bibitem{b12}
I.~Protasov, S.~Slobodianiuk,
{\it Thin subsets of groups}, Ukr. Math. J., submitted.

\bibitem{b13}
T.~Tao, V.~Vu,
{\it Additive Combinatorics}, Cambridge University Press, 2006.

\end{thebibliography}
\end{document}